\documentclass[11pt]{article}
\usepackage{latexsym,amsmath,amssymb,theorem}
 \usepackage[applemac]{inputenc}  
\usepackage[all]{xy}

\usepackage{color} 
\usepackage{hyperref}
\usepackage{mathbbol}
\newenvironment{Proof}{\noindent \bf Proof   \rm}{\hspace*{\fill}
$\square$ \medskip}

{\theorembodyfont{\itshape}
\newtheorem{Lemma}{Lemma}
\newtheorem{Proposition}[Lemma]{Proposition}
\newtheorem{Fact}[Lemma]{Fact}
\newtheorem{Facts}[Lemma]{Facts}

\newtheorem{Corollary}[Lemma]{Corollary}
}
{\theorembodyfont{\upshape}

\newtheorem{Definition}[Lemma]{Definition}
\newtheorem{Remark}[Lemma]{Remark}

}

\newcommand{\ot}{\otimes}

\newcommand{\N}{\mathbb{N}}
\newcommand{\Z}{\mathbb{Z}}

\newcommand{\BC}{\mathcal{C}}

\newcommand{\BBT}{\mathbb{T}}

\newcommand{\BA}{\mathcal{A}}

\newcommand{\Cat}{\mathit{Cat}}

\newcommand{\op}{\mathsf{op}}
\newcommand{\ob}{\mathsf{ob}}

\newcommand{\id}{\mathsf{id}}

\newcommand{\sm}{\mathsf{M}}

\newcommand{\Set}{\mathit{Set}}

\newcommand{\Mon}{\mathit{Mon}}

\newcommand{\lra}{\rightarrow}

\newcommand{\wra}{\xymatrix{ \ar@{~>}[r] & }}

\newcommand{\xra}{\xrightarrow}

\newcommand{\spant}[6]{
$\xymatrix@=1,5em{& 
{#1 }\ar[dl]_{#4 } \ar[dr]^{#5} &\\
 {#2}\ar[rr]_{#6 } & &{#3} }$
}
\newcommand{\spanl}[5]{$ #2   \xleftarrow{#4} #1\xrightarrow{#5} #3$}

\newcommand{\spano}[5]{
$\xymatrix@=1em{& 
{#1 }\ar[dl]_{#4 } \ar[dr]^{#5} &\\
 {#2} & &{#3} }$
}

\makeatletter
\newcommand{\xRightarrow}[2][]{\ext@arrow 0359\Rightarrowfill@{#1}{#2}}
\makeatother

\begin{document}

\title{The algebra of Feistel-Toffoli schemes}

  \author{ Laurent Poinsot\\\small CREA, \'Ecole de l'Air et de l'Espace, Salon-de-Provence and\\
\small LIPN, UMR CNRS 7030, {Sorbonne Paris North University, Villetaneuse, France. }\\ 
{\tt \small laurent.poinsot@lipn.univ-paris13.fr}
 \and  Hans-E. Porst\thanks{Permanent address: Besselstr. 65, 28203  Bremen, Germany.} \\
\small Department of Mathematical Sciences, University of Stellenbosch, \\ \small Stellenbosch, South Africa. \\
 {\tt\small porst@uni-bremen.de} 
 }

\date{}

\maketitle

\begin{abstract} 
The process of replacing an arbitrary  Boolean function by a bijective one, a fundamental tool in reversible  computing and in cryptography, is interpreted algebraically as a {particular instance of a certain} group homomorphism from the $X$-fold cartesian power of a group $G$ into the automorphism group of the free $G$-set over the set $X$.  It is shown that this construction {not only can be generalized from groups to monoids but, more generally,} to internal categories in arbitrary  finitely complete categories where it becomes a cartesian isomorphism between certain discrete fibrations.
\end{abstract}

\section{Introduction}

 The basic idea behind the so-called  Toffoli gate, a fundamental tool in reversible computing and, hence, one of the ingredients to prove that every classical program may be executed by a quantum computer, and called the \em Fundamental Theorem \em in \cite{Toffoli} is the following fact: 
 every function $f\colon 2^m\to 2^n$  can be  turned into some bijective function $E_m
 (f) \colon 2^m\times 2^n = 2^{m+n}\to 2^{m+n}= 2^m\times 2^n$ 
 in such a way that $f$ may be retrieved from it.
 Thinking of  $2 =\{0,1\}$ as $\Z_2$ and using the notations  $\vec{x}:=(x_1,\cdots,x_m)\in \mathbb{Z}_2^m$ and $f(\vec{x})=:(f_1(\vec{x}),\cdots,f_n(\vec{x}))$ the map $E_m(f)\colon \mathbb{Z}_2^{m+n}\to\mathbb{Z}_2^{m+n}$ is  given by $$(\vec{x},\vec{y})\mapsto(\vec{x},f_1(\vec{x})+y_1,\cdots,f_n(\vec{x})+y_n).$$ 
 About a decade older is the consideration of its special instance with $n=m$; 
this is of  interest in   cryptography, since it provides a method for developing the general architecture of some symmetric encryption schemes commonly known as {\em Feistel ciphers}~\cite[Definition~7.81]{Menezes}.

Hence, the natural question arises  whether this map $E_m$ has any conceptual algebraic meaning and, if so, whether this might make sense in a more general context.
In order to    {find an answer to this question} let us observe first that in the definition of $E_m(f)$ neither the    {algebraic}  structure of the domain $\mathbb{Z}_2^m$ of $f$ is used nor the fact that the codomain of $f$, that is, the additive monoid $\mathbb{Z}_2^n$ is a group. This leads us to start simply with a map $X\xra{f}M$ where $X$ is a set and $\mathsf{M}=(M,\cdot,1)$ is a monoid such that we obtain    {by the construction above} a map  $X\times M\xra{E_X(f)}X\times M$, $(x,m)\mapsto (x,f(x)\cdot m)$. Noting that $X\times M$ is the (underlying set of the) free  object $F_\sm X$ over $X$ in the category $\Set_\sm$ of right $M$-acts   one immediately sees that the map $E_X(f)$ is nothing but the homomorphic extension $\langle id_X, f\rangle^\sharp$ of the map $X\xra{\langle id_X, f\rangle}  X\times M  $ with $x\mapsto (x,f(x))$, that is, the unique map making the following diagram commute, where $\eta_X(x) = (x,1)$.
\begin{equation*}\label{diag:1}
\begin{aligned}
\xymatrix@=2.3em{
     X\ar[r]^{\eta_X  }\ar[dr]_{ \langle id_X, f\rangle }&{X\times M } \ar[d]^{E_X(f)=\langle id_X, f\rangle^\sharp }  \\
         &   X\times M 
}
\end{aligned}
\end{equation*}
This not only answers the question at the beginning of   the previous paragraph but also 
 describes the reconstruction of $f$ from $E_X(f)$ structurally:  
 \begin{equation}\label{eqn:1}
f= X\xra{\eta_X}X\times M\xra{E_X(f)}X\times M\xra{\pi_M}M =:L_X(E_X(f)).
\end{equation}
Considering the map $f\colon X\to M$  as an element of the $X$-fold cartesian power of $M$  the map $E_X$ becomes a map $M^X\xra{}\Set_\sm(F_\sm X,F_\sm X)$. Denoting by $s\ast t$  the multiplication in $\sm^X$ one    {gets by a straightforward calculation the equations 
$\langle id_X,s\ast t\rangle^\sharp = \langle id_X, s \rangle^\sharp\circ \langle id_X, t\rangle^\sharp \text{\ and \ } \langle id_X,1_{\sm^X}\rangle^\sharp = id_{F_\sm X}$, for any pair $(s,t)\in M^X$}. Thus,  
$E_X$ even is a monoid homomorphism from the  $X$-fold power $\sm^X$ of  $\sm$ in the category $Mon$ of (small) monoids  into the endomorphism  monoid $\Set_\sm(F_\sm X,F_\sm X)$.  
{Since $E_X$ is injective by Equation \eqref{eqn:1} one concludes that $\sm^X$ is isomorphic to a submonoid of this endomorphism monoid. }

In view of the starting examples, where $   {E_X(f)}$ was bijective, we note that this homomorphism property  implies that the monoid $E_X(\sm^X)$ is a subgroup of the automorphism group  of the free $\sm$-act over $X$, provided that the monoid $\sm$  and, hence, the monoid $\sm^X$  even is a group.

This construction will be shown to be  functorial in the following sense: 
the assignments $X\mapsto \sm^X$ and $X\mapsto E_X(\sm^X) $ define contravariant functors $P^\ast_\sm$ and $Q^\ast_\sm$  from $\Set$ to $\Set$ which are naturally isomorphic by means of the family $(E_X)_X$. 
This is in other words: the natural transformation $(E_X)_X$ determines a cartesian isomorphism between  the discrete fibrations $P_\sm$ and $Q_\sm$ over $\Set$ whose base change functors 
are the functors $P^\ast_\sm$ and $Q^\ast_\sm$. 

All of this generalizes quite obviously to any category $\BC$ with finite products instead of the category $\Set$ of (small) sets  as will be shown in Section \ref{sec:2}. Having in mind that a monoid is a one-object internal category in $\Set$ one may wonder whether even a generalization to \em arbitrary \em internal categories in \em arbitrary \em finitely complete categories $\BC$ is possible; and this the more so, since an internal category in  $\BC$, which has $O$ as its object of objects,   is the same as a monoid in a certain monoidal category  $\BC_O$; since $\BC_O$ contains an isomorphic copy of $\BC/O$ as a full monoidal subcategory, one  is lead to the following conjecture.\footnote{ $\BC_O$ is the monoidal category $\BC_O:=Span(\BC)(O,O)$, where $Span(\BC)$ denotes the bicategory of spans in $\BC$ (see below).  {Note that the categories $Span(\Set)(1,1)= Set_1$, $\Set/1$ and $\Set$ are isomorphic. }} \em All of the above can be generalized to pairs $(X,\sm)$ where $\sm$ is  an internal category in  {a finitely complete category} $\BC$  and $X$ is an object of $\BC/O$. \em  The proof of this conjecture then is the purpose of this note, where particular emphasis is given to the question in which sense the occurring constructions are functorial. The main part of this note, Section \ref{sec:3},  contains the proof of this conjecture, where particular emphasis is given to the question in which sense the occurring constructions are functorial.

\section{The case of internal monoids}\label{sec:2}

Let  $\BC$ be a cartesian category, that is, a category  with finite products and  $\sf{M} = $ $(M,M\times M\xra{m}M, 1\xra{e}M)$  a monoid in $\BC$. Then the following hold for any $\BC$-object $C$:
\begin{enumerate}
\item 
$(\BC(C,M), \ast, e_C)$ is a monoid, for each $\BC$-object $C$ with 
$e_C:= C\xra{!}1\xra{e}{M}$ where $1$ is the terminal object of $\BC$, and $\phi\ast\psi := C\xra{\Delta_C}C\times C\xra{\phi\times\psi} {M}\times {M}\xra{m} {M} $. This is a special instance of the so-called \em convolution monoid \em given by a comonoid $C$ and a monoid $M$ in a monoidal category (see also Section \ref{sec:Mval}), when considering $(\BC,\times, 1)$ as a cartesian monoidal category, that is, a monoidal category whose monoidal structure is given by binary products and a terminal object.
 Note, that in case $\BC =\Set$, this monoid is simply the $C$-fold power of the monoid $\sm$.
 \\
 If $\sf{M} $ = $ (M,M\times M\xra{m}M, 1\xra{e}M, M\xra{i}M)$  even is a group in $\BC$ then  $(\BC(C,M), e_C, \ast)$ is a group, since $\phi\ast(i\circ\phi) = e_C = (i\circ\phi)\ast\phi$, for each $\phi\in \BC(C,M)$.
\\ 
  For a different argument using the Yoneda functor see  e.g. \cite[III.6]{MacL}.
\item There is the category $\BC_\sm$ of $\sf{M}$-acts $(C,C\times M\xra{\alpha}C)$ in $\BC$ defined in  the obvious way and an obvious forgetful functor $|-|\colon\BC_\sm\xra{}\BC$. This again is a special instance of the category of right $\sm$-modules with respect to a monoid in a monoidal category.
\item The forgetful functor $|-|$ of the category  of right $\sm$-modules  in a monoidal category $\BC$ has a left adjoint
  $F_\sm$ (see e.g. \cite{Par}). This specializes to a category with finite products as follows.
\begin{enumerate}
\item 
$|-|\colon\BC_\sm\xra{}\BC$ has a left adjoint $F_\sm$ given by $$F_\sm(C) = C\times M\times M\xra{\id_C\times m}C\times M $$
 with units  $C\xra{\eta_C}C\times M    = C\simeq C\times 1 \xra{id_C\times e} |(C\times M, \alpha)| $. 
\item The underlying  $\BC$-morphism of the homomorphic extension $f^\sharp$ of a $\BC$-morphism $C\xra{f} Y\times M$ is given by  
$ C\times M \xra{f\times id_M}Y\times M\times M\xra{ id_Y\times m}Y\times M$.
\end{enumerate}
Note that Items (a) and (b), when specialized to $\Set$, are precisely the descriptions of the units and homomorphic extensions in the introduction.
\end{enumerate}
Consequently, the arguments used to prove the properties of the map $E_X$ in the introduction apply literally  in the case where the category under consideration is not $\Set$ but an arbitrary category with finite products. We so obtain the 
\begin{Proposition}\label{prop:}
Let $\sf{M} = $ $(M,M\times M\xra{m}M, 1\xra{e}M)$  be a monoid in   a  category $\BC$ with finite products. Then the assignment $(C\xra{f}M)\mapsto   \langle id_C, f\rangle^\sharp$
defines a  map $E_C\colon \BC(C,M)\xra{}\BC_\sm(F_\sm C, F_\sm C)$ such that the following hold.
\begin{enumerate}
\item $E_C$ is monoid morphism from the convolution monoid into the $\BC_\sm$-endomorphism monoid.
\item $E_C$ is  a section.
\item If $\sm$ even is a group in $\BC$ then, for each $C\xra{f}M$, the endomorphism $E_C(f)$ is an automorphism.
\end{enumerate}
\end{Proposition}
\begin{Remark}\label{rem:}\rm
One certainly could have shown the above without explicitly referring to the theory of monoidal categories. We prefer our approach for the following reasons: (a) the monoidal methods will have to be used below in any case.  (b) When interpreting varieties (in the sense of universal algebra) as those of monoids, groups and $\sm$-acts in a category with finite products one cannot  prove the existence of free algebras in general. Hence, the existence of free $\sm$-acts used above is rather untypical in this context   --- it becomes natural, however, when considering a category with finite products as a cartesian \em monoidal \em category.
\end{Remark}

\section{The case of internal categories}\label{sec:3}
As from now the category $\BC$ is supposed to be finitely complete.

\subsection{Prerequisites}

\subsubsection{Some notation}\label{sec:not}
 Given morphisms $M\xra{c}O$ and $N\xra{d}O$ in $\BC$, the  pullback of $c$ along $d$ over $O$  will be denoted by
 \begin{equation*}\label{st}
\begin{aligned}
\xymatrix@=2em{M\!{\, _d\times_c }\, N  \ar[r]^{\ \ \pi^M_c}\ar[d]_{\pi^N_d} & N\ar[d]^c   \\
M  \ar[r]_{ d}& O
}
\end{aligned}
\end{equation*}
If possible, we will write only $\pi_c$ and $\pi_d$  instead of  $\pi^M_c$ and $\pi^N_d$, respectively; occasionally we may simply write $p$ instead of $\pi_c$ and $q$ instead of $\pi_d$.

Note that by the standing assumption on $\BC$ the object $M{_f\times_g }N$ is a (regular) subobject of $M\times N$.
Given a pair of morphisms $M'\xra{f }M$  and $N'\xra{g }N$ such that  the morphism $M'\times N'\xra{ f\times g} M\times N$ factors over $M\!{\, _d\times_c }N$ we will  by slight abuse of notation its corestriction $M'{\times }N'\lra M{\,_ d\times_c }N$ simply denote  by $f\times g$ as well;  {more generally,  any restriction of this morphism to a subobject $X$ of $M'\times N'$ which factors over $M{\,_ d\times_c }N$ may be denoted by $X\xra{f\times g} M{\,_ d\times_c}N$}. Similarly, if $N'=M'$ and the morphism $M'\xra{ \langle f,g\rangle} M\times N$ factors over  $M{\,_ d\times_c }N$, that is, if $d\circ f =c\circ g$,   the resulting morphism $M'\xra{ } M{\,_ d\times_c}N$ will be denoted simply by $ \langle f,g\rangle$ as well.

\subsubsection{The monoidal category $\BC_O$ }
The categories described  in this section, to which we referred already in the introduction, are of crucial importance when working with internal categories.

Let $O$ be a fixed object in the finitely complete category $\BC$. Then the category $\BC_O$ has as its objects all pairs  $(x\xra{ f}O, x\xra{ g}O)$ of morphisms in  $\BC$, called \em spans \em and usually denoted by \spanl{x}{O}{O}{f}{g}, and as its morphisms, also referred to as {\em $2$-cells} (see Remark~\ref{rem:span} below),  \spanl{x}{O}{O}{f}{g} $\Rightarrow$ \spanl{x'}{O}{O}{f'}{g'}  all $\BC$-morphisms  $x\xra{t}x'$ making the following diagram commute.
\begin{equation*}\label{diag:2}
\begin{aligned}
\xymatrix@=1em{& x \ar[dr]^{ g }\ar[dl]_{ f }\ar[dd]_{ t } &\\
    O &{ }&   O\\
 &x' \ar[ur]_{ g'} \ar[ul]^{ f' }         &     
}
\end{aligned}
\end{equation*}
The composition in $\BC_O$ is the composition in $\BC$ and the identity on  \spanl{x}{O}{O}{f}{g} is the identity $id_x$ in $\BC$. The functor $\BC_O\xrightarrow{|-|}\BC$ denotes the forgetful functor given by the assignment $((O\xleftarrow{f}x\xrightarrow{g}O)\xrightarrow{t}(O\xleftarrow{f'}x'\xrightarrow{g'}O))\mapsto (x\xrightarrow{t}x')$. Note that when $\BC$ is locally small then $\BC_O$ is locally small as well. The category $\BC_O$ is equipped with a monoidal structure defined as follows. 
\begin{enumerate}
\item $\text{(\spanl{y}{O}{O}{h}{k})} \ot \text{(\spanl{x}{O}{O}{f}{g})} := \text{\spanl{x {_{g}\times_{h}} y }{O}{O}{f\circ \pi_g}{k\circ \pi_h}}$ is given by pullbacks as visualized by the following diagram, where we occasionally denote the $\BC$-object $x {_g\times_h}y$ by $x{\times_O}y$.\\ 
{\small 
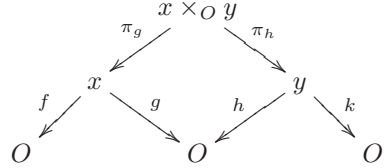
$$\xymatrix@=1.5em{& & x\times_O y\ar[dl]_{\pi_g} \ar[dr]^{\pi_h}& & \\
& x \ar[dr]^{ g }\ar[dl]_{ f } &  & y\ar[dr]^{ k }\ar[dl]_{ h } &\\
    O &{ }&   O & & O
}$$}
\item Given $\BC_O$-morphisms
$s\colon$ \spanl{y}{O}{O}{k}{l}  $\Rightarrow$ \spanl{y'}{O}{O}{k'}{l'}  and $t\colon $ \spanl{x}{O}{O}{f}{g} $\Rightarrow$ \spanl{x'}{O}{O}{f'}{g'}  then $t\ot s$ is the $\BC$-morphism from $x\times_O y$  to $x'\times_O y'$ induced by the pullback property of $x'\times_O y'$ from the $\BC$-morphisms $x\times_O y \xra{ p}x\xra{ t} x' $ and $x\times_O y \xra{ q}y\xra{ s} y'$ as visualized by the commutative diagram below. In other words, $t\ot s$ is the unique morphism with $p'\circ (t\ot s) = t\circ p$ and 
$q'\circ  (t\ot s) = s\circ q$. 
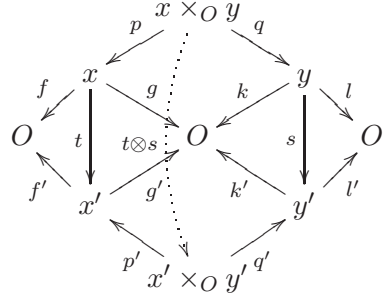
\begin{equation*}\label{diag:ast}
\begin{aligned}
\xymatrix@=1em{& & x\times_O y\ar[dl]_{p} \ar[dr]^{q}\ar@/_1pc/@{.>}[dddd]_{t\ot s}& & \\
& x \ar[dr]^{ g }\ar[dl]_{ f }\ar[dd]_{ t } &  & y\ar[dr]^{ l }\ar[dl]_{ k }\ar[dd]_{ s } &\\
    O &{ }&   O & & O\\
 &x' \ar[ur]_{ g'} \ar[ul]^{ f' }         &     &y' \ar[ur]_{ l'} \ar[ul]^{ k' }  &\\
 & & x'\times_O y' \ar[ul]^{p'} \ar[ur]_{q'}& &
}
\end{aligned}
\end{equation*}
\item The monoidal unit is the span $\bar{O}$: = \spanl{O}{O}{O}{id}{id}.
\end{enumerate}
The associativity and units constraints are easily obtained from the pullback property.
\begin{Remark}\label{rem:span}\rm
The fact above is a consequence of the well known results due to B\'enabou \cite{Ben}, that the  spans  of a finitely complete category $\BC$ form a bicategory, and that the categories of 1-cells of  a bicategory are monoidal (see also \cite{DPP}).
\end{Remark}

\subsubsection*{The cartesian subcategory ${\widehat{\BC}_O}$}
In the sequel we denote $\BC_O$-objects of the form  \spanl{A}{O}{O}{f}{f} simply by  $\mathsf{f}_A$ and by ${\widehat{\BC}_O}$ the full subcategory of $\BC_O$ spanned by all the spans 
$\mathsf{f}_A $. This category is obviously isomorphic  to the slice category $\BC/O$ by means of the functor $D\colon\BC/O\xra{}\BC_O$ with $  \big((X\xra{f} O)\xra{\phi}(Y\xra{g}O)\big)\mapsto  (\mathsf{f}_X\xra{\phi}\mathsf{g}_Y)$. We, hence, may not distinguish notationally between (objects of) $\BC/O$ and 
${\widehat{\BC}_O}$.
 
 For ${\widehat{\BC}_O}$-objects $\mathsf{f}_A $ and $\mathsf{g}_B $ one has 
  $\mathsf{f}_A \ot \mathsf{g}_B = O \xleftarrow{f}A\xleftarrow{\pi_f} A{_f\times_g} B\xra{\pi_g}B\xra{g}O$,  
 where $f\circ\pi_f = g\circ\pi_g$. Hence, the $\BC_O$-object  $\mathsf{f}_A \ot \mathsf{g}_B = (\mathsf{f\circ\pi_f})_{A{_f\times_g} B}$  belongs to ${\widehat{\BC}_O}$. 
 Since the $\BC$-morphisms $(A{_f\times_g} B\xra{\pi_f}A)$ and $(A{_f\times_g} B\xra{\pi_g}B)$ are the product projections of $(A\xra{f} O)\times (B\xra{g}O)$ in $\BC/O$ one obtains
 
  \begin{Lemma}\label{lem:crux}
${\widehat{\BC}_O}$ is a full monoidal subcategory of $\BC_O$ and $D\colon\BC/O\xra{}\BC_O$ is a monoidal equivalence between the cartesian category $\BC/O$ and the full monoidal subcategory ${\widehat{\BC}_O}$ of $\BC_O$. 
In more detail:
given  ${\widehat{\BC}_O}$-objects $\mathsf{f}_A $ and $\mathsf{g}_B $ the following hold.
\begin{enumerate}
\item $(\mathsf{f}_A \ot \mathsf{g}_B,\pi_f,\pi_g) $ is  a product of $\mathsf{f}_A $ and $\mathsf{g}_B $ in ${\widehat{\BC}_O}$. 
\item Since  for each  $\mathsf{f}_A$ = \spanl{A}{O}{O}{f}{f} in ${\widehat{\BC}_O}$ the $\BC$-morphism $f$ is the unique ${\widehat{\BC}_O}$-morphism $\mathsf{f}_A\xra{ !} $
(\,\spanl{O}{O}{O}{id}{id}), the ${\widehat{\BC}_O}$-object  $\bar{O} = \mathsf{id}_O$ is terminal in  ${\widehat{\BC}_O}$. 
\end{enumerate}
\end{Lemma}

 \subsubsection*{Some calculation rules}\label{fact:pseudo}

Given spans  $\bar{X} =$ (\spanl{X}{O}{O}{x}{y}) and  $\bar{M} =$ (\spanl{M}{O}{O}{d}{c}) there is, for any pair of $\BC$-morphisms  $A\xra{ \alpha}M$ and $A\xra{\xi }X$  satisfying the condition $d\circ\alpha = y\circ\xi$,  
 the unique $\BC$-morphism $A\xra{ \langle \xi,\alpha\rangle} |\bar{X}\ot\bar{M}|$ with $\pi_d\circ \langle \xi,\alpha\rangle =\alpha$ and $\pi_y\circ \langle \xi,\alpha\rangle=\xi$   (see Section \ref{sec:not}). 
 The following facts concerning the interplay of the monoidal structure of $\BC_O$ and the cartesian structure of its subcategory ${\widehat{\BC}_O}$ are easy to prove and will be of use later.

\begin{enumerate}
\item If  $\alpha$ and $\xi$  are 2-cells $\mathsf{f}_A\Rightarrow  \bar{M}$  and $\mathsf{f}_A\Rightarrow\bar{X}$, respectively, then $\langle \xi,\alpha\rangle$ exists and is a 2-cell $\mathsf{f}_A\Rightarrow  \bar{X}\ot \bar{M}$. 

\item Denoting for $\mathsf{f}_A$ in ${\widehat{\BC}_O}$ by $\mathsf{f}_A\xra{\Delta}\mathsf{f}_A\ot\mathsf{f}_A$ or $\Delta_A$, if necessary, its diagonal (w.r.t. the cartesian structure of ${\widehat{\BC}_O}$), for any pair of $\BC_O$-morphisms $\mathsf{f}_A\xra{\xi}\bar{X}, \mathsf{f}_A\xra{\alpha}\bar{M} $   one has $\langle \xi, \alpha\rangle = (\xi\ot \alpha)\circ\Delta $. 
\item For any pair  of $\BC_O$-morphisms $\alpha,\beta\colon \mathsf{f}_A\xra{}\bar{M}$  the following diagram commutes. 
\begin{equation*}\label{diag:3}
\begin{aligned}
\xymatrix@=2.5em{
 \mathsf{f}_A \ar[r]^{\Delta}\ar[d]_{\Delta}&  \mathsf{f}_A\ot \mathsf{f}_A  \ar[r]^{\Delta\ot \beta  }&{ \mathsf{f}_A\ot \mathsf{f}_A\ot \bar{M}} \ar[d]^{ id\ot\alpha\ot id }  \\
\mathsf{f}_A\ot \mathsf{f}_A\ar[r]_{id\ot \Delta } &  \mathsf{f}_A\ot \mathsf{f}_A\ot \mathsf{f}_A   \ar[r]_{  id\ot \alpha\ot\beta}          & \mathsf{f}_A\ot  \bar{M}\ot \bar{M}
}
\end{aligned}
\end{equation*}

\item Concerning the construction $\mathsf{f}_A\ot \bar{M}$ we note moreover: 
\begin{enumerate}
\item the pullback projection $\mathsf{f}_A\ot \bar{M}\xra{\pi^A_d }\bar{M}$ is a 2-cell; 
\item for each 2-cell $\alpha\colon \mathsf{f}_A\Rightarrow  \mathsf{f}_A\ot \bar{M}$ the composite $\pi_f\circ\alpha$ is a 2-cell $\mathsf{f}_A\Rightarrow\mathsf{f}_A$ ({though the pullback projection $\pi_f$ will fail to be a 2-cell in general)};
\item for any 2-cell $\phi\colon \mathsf{f}_A\Rightarrow \mathsf{g}_B$ the following diagrams commute
\begin{equation*}\label{7}
\begin{aligned}
\xymatrix@=2.5em{
 \mathsf{f}_A\ot \bar{M}\ar[r]^{\pi^A_d}\ar[d]_{\phi\ot id} &\bar{M}\ar[d]^{id} \\
\mathsf{g}_B\ot \bar{M}\ar[r]^{\pi^B_d}  &\bar{M}
}
\qquad \hspace{1cm}
\xymatrix@=2.5em{
  |\ \mathsf{f}_A\ot \bar{M}|\ar[d]_{|\phi\ot id|}  \ar[r]^{\pi_f} &|\mathsf{f}_A|
  \ar[d]^{|\phi|}  \\ 
 |\mathsf{g}_B\ot \bar{M}| \ar[r]_{\pi_g} &|\mathsf{g}_B|
   }
\end{aligned}
\end{equation*}
\end{enumerate}
\end{enumerate}

 \subsubsection{Internal categories and groupoids}\label{sec:intcat}
 \subsubsection*{Internal categories}\label{intcat}
 
Recall  (see e.g. \cite[Ex. 3C]{AHS}) that a small category {can be seen as} a sixtuple $${\mathsf{M}}=(O,M,d,c\colon M\xra{ }O, O\xra{\eta }M,M\!\times_O\! M\xra{\mu }M)$$ 
 of sets and maps, where $O$ and $M$ are   the sets of objects and morphisms, respectively, $c$ and $d$ are the  (co)domain maps,  ${\eta}$  is thought of as the  family of identities, and $\mu$  is the composition map. These data must satisfy the category axioms, that is, they make the obvious  diagrams commute. 
This all makes sense if the terms \em sets \em and  \em maps \em are replaced by \em objects \em and \em morphisms \em of a finitely complete category $\BC$. In this case then $\mathsf{M}$ is called an \em internal category in $\BC$\em. Obviously an internal category with $O=1$, the terminal object of $\BC$, is essentially the same thing as an internal monoid $(M, \mu,\eta)$ in $\BC$. 

The following generalizes a familiar fact about internal monoids (see e.g. \cite[III.6]{MacL}) and is 
 easy to see. The Yoneda functor $Y\colon\BC\xra{}\Set^{\BC^\op}$ maps an internal  category $\mathsf{M} =(O,M,d,c,\mu,\eta)$ in $\BC$ (here assumed locally small) to an internal category in $\Set^{\BC^\op}$ and each evaluation functor $\Set^{\BC^\op}\xra{ev_C}\Set$ maps this to the small  category $\mathsf{M}_C$ which has  
 \begin{itemize}
\item $\BC(C,O)$ as its set of objects and $\BC(C,M)$ as its set of morphisms,  
\item  the maps ${d}_C:=\BC(C,M)\xra{\BC(C,d)}\BC(C,O)$ and ${c}_C:= \BC(C,M)\xra{\BC(C,c)}\BC(C,O)$ as domain and codomain assignments and the map $\BC(C,O)\xra{\BC(C,\eta )}\BC(C,M)$ as its map assigning units.
\item the map $\BC(C,M)\times_{\BC(C,O)} \BC(C,M)= \BC(C,M\times_O M)\xra{\BC(C, \mu)}\BC(C,M)$ as its composition, where $\BC(C,M)\times_{\BC(C,O)} \BC(C,M)$ is  the pullback of ${c_C}$ along ${d}_C$ in $\Set$.
\end{itemize}
In other words, the category $\mathsf{M}_C$ has
 \begin{itemize}
 \item the $\BC$-morphisms $C\xra{ f}O$ as its objects, 
\item the sets $\{ C\xra{\alpha} M \mid {d}_C (\alpha) = d\circ\alpha = f \text{\ and \ } {c}_C(\alpha) =  c\circ\alpha=g\}$ as its hom-sets $\hom((C,f),(C,g))$,
\item the $\BC$-morphisms $id_{(C,f)}  := C \xra{f }O\xra{\eta }M $ as its units,
\item the composition $((C,g)\xra{ \beta}(C,h))\circ ((C,f)\xra{ \alpha}(C,g))= C\xra{\langle\beta ,\alpha\rangle} M\times_OM \xra{\mu }M$.
\end{itemize}
  
\begin{Remark}\label{rem:catfibr} 
Given an internal category $\mathsf{M} = (O,M,d,c,\mu,\eta)$ in $\BC$, the assignment $C\mapsto \sm_C$ can obviously be extended to a functor $\BC^\op\xra{}\Cat$ (where $\Cat$ denotes the category of (small) categories); consequently, there exists a split fibration ${U}\colon \BC_{(O,M)}\xra{}\BC$    {called the \em externalization of $\mathsf{M}$} (see e.g. \cite{Jacobs}), whose fibres then are the categories $\mathsf{M}_C$ just defined.
\end{Remark}

Recall that an \em internal functor \em  $(O,M,c,d,\eta,\mu) \lra (O',M',c',d',\eta',\mu')$ between internal categories in $\BC$  is a pair of $\BC$-morphisms $F = (O\xra{ F_o}O', M\xra{ F_m}M')$ such that the following diagrams commute.
\begin{equation*}\label{}
\begin{aligned}
\xymatrix@=2.5em{ M  \ar@<.5ex>[r]^{c }\ar@<-0.5ex>[r]_{d }\ar[d]_{ F_m} & O\ar[d]_{F_o}\ar[r]^{\eta}  & M\ar[d]^{F_m}  \\
 M'   \ar@<.5ex>[r]^{c' }\ar@<-0.5ex>[r]_{d' }& O'\ar[r]^{\eta'} & M'
}
\qquad 
\xymatrix@=2.5em{ M\times_OM  \ar[d]_{ F_m\times F_m  }   \ar[r]^{\ \ \mu} & M \ar[d]^{F_m} \\
 M'\times_{O'}M'\ar[r]^{\ \ \mu' }&M'}
\end{aligned}
\end{equation*}
One so obtains ${Cat}_{in}(\BC)$, the \em category of internal categories and internal functors in $\BC$\em. 
\begin{Remark}\label{rem:Catin}\rm
Obviously, every finite limit preserving  functor $\mathcal{C}\xra{K}\mathcal{C}'$ between finitely complete categories  induces a functor $Cat_{in}(K)\colon {Cat}_{in}(\mathcal{C})\to {Cat}_{in}(\mathcal{C}')$. 
\end{Remark}

\subsubsection*{Internal categories  in the language of spans }
In the language of spans the description of an internal category in a finitely complete category $\BC$ reads as follows: an internal category  $\mathsf{M}=(O,M,d,c,\eta, \mu) $   is a triple $(\bar{M},\eta,\mu)$ where 
\begin{itemize}
\item  $\bar{M}$ = \spanl{M}{O}{O}{d}{c}   is  an object of $\BC_O$, 
\item   $\eta\colon  \bar{O}$ = \spanl{O}{O}{O}{id}{id} $\Rightarrow \bar{M}$  and  $\mu\colon \bar{M}\ot \bar{M} \Rightarrow \bar{M}$ 
are morphisms in $\BC_O$,
 \item the following diagrams commute 
\begin{equation*}\label{ax}
\begin{aligned}
\xymatrix@=2em{\bar{M}\ot\bar{M}\ot\bar{M} \ar[rr]^{\ \ \ \ \mu\ot id }\ar[d]_{id\ot\mu } && \bar{M}\ot\bar{M}\ar[d]^{\mu}    \\
 \bar{M}\ot \bar{M} \ar[rr]_{\mu  }&& \bar{M}
}
\qquad 
\xymatrix@=2em{ \bar{O}\ot\bar{M} \ar[r]^{\eta\ot id }\ar[r]^{}\ar[dr]_{\pi_d }&\bar{M}\ot\bar{M}\ar[d]^{\mu} &  \bar{M}\ot\bar{O}\ar[dl]^{\pi_c}\ar[l]_{id\ot \eta}   \\
   & \bar{M}
}
\end{aligned}
\end{equation*}
\end{itemize}
 Since these equations  are nothing but the axioms for the triple $(\bar{M},\mu,\eta)$ with $\bar{M} =$ (\spanl{M}{O}{O}{d}{c}) to be a  monoid in the monoidal category $\BC_O$, we have got the following fact.

\begin{Fact}\label{fact:int}\rm
If $\BC$ is a finitely complete category  then  internal categories  in $\BC$ with $O$ its object of objects are precisely the monoids in the monoidal category   $\BC_O$.\footnote{{Alternatively one may say that these internal categories are monads in the bicategory $Span(\BC)$ on the object $O$ in the sense of \cite{Street}; this view however is not the appropriate one for our purpose.}}
\end{Fact}

\subsubsection*{Internal groupoids}
As groups and groupoids can equivalently be considered as  monoids and categories, respectively, with an additional property (every element has an inverse) or an additional structure (the map mapping an element to its inverse), the same holds for internal groupoids. Here, however, neither the additional property nor the respective equivalence is as obvious as in the cases just mentioned (see e.g. \cite{Illusie}).
We, hence, refrain from describing internal groupoids by a particular property of an internal category and only  use the following classical definition (see \cite{Brown}) which suits our needs best. 

\begin{Definition}\label{def:str}\rm
Let $\BC$ be a finitely complete category. An   \em internal groupoid in $\BC$ \em is a pair $(\mathsf{M},\iota)$ where $\mathsf{M}= (O,M,d,c,\eta,\mu)$ is an internal category in $\BC$  and  $\iota\colon M\xra{ }M$ is 
a $\BC$-morphism such that 
\begin{enumerate}
\item $c\circ \iota = d$ and $d\circ \iota = c$ and
\item $\mu\circ\langle id,\iota\rangle =\eta\circ c$  and $\mu\circ {\langle \iota, id\rangle} = \eta\circ d$.
\end{enumerate}
\end{Definition}

If now $(\mathsf{M},\iota) = (O,M,d,c,\eta,\mu,\iota)$ is an internal groupoid then, for each $\BC$-object $C$, 
the category  $\mathsf{M}_C$ is equipped with  a map ${\iota}_C =\BC(C,\iota)\colon \BC(C,M)\xra{ }\BC(C,M) $
{satisfying the conditions
\begin{enumerate}
\item ${c_C}\circ {\iota}_C  = {d_C}$ and  ${d_C}\circ {\iota}_C  = {c_C}$ and 
\item the following diagram commutes
\begin{equation*}\label{diag:4}
\begin{aligned}
\xymatrix@=2.5em{
{ \BC(C,M)}\ar[r]^{ \BC(C,\eta\circ c)} \ar[dr]_{\BC(C,\langle id,\iota\rangle)}&\BC(C,M) & \BC(C,M) \ar[dl]^{\BC(C,\langle \iota, id\rangle) }\ar[l]_{\BC(C,\eta\circ d)}     \\
& {\BC(C,M\times_OM)} \ar[u]_{\BC(C, \mu)}    
}
\end{aligned}
\end{equation*}
\end{enumerate}
We, hence, can conclude that every morphism $\alpha$ in $\mathsf{M}_C$ is an isomorphism with inverse $\alpha^{-1} = \iota_C(\alpha)$, such that $\mathsf{M}_C$ in fact is a groupoid and that ${{\iota}_C}\circ{\iota}_C=id_{\BC(C,M)}$, for each $\BC$-object $C$. Since the family $({{\iota}_C})_C$ is the image of the full and faithful Yoneda embedding $Y$ the $\BC$-morphism $\iota$ is an idempotent isomorphism as well. 
We so have got the following results.
\begin{Facts}\label{prop:onetoone}\rm 
If  $(\mathsf{M},\iota)$ is an internal groupoid in $\BC$ then all categories $\sm_C$ are groupoids and $\iota\circ\iota =id_M$.
\end{Facts}

\subsubsection{Discrete fibrations}\label{discfib}
We recall the following well known concepts and facts (see e.g. \cite{Riehl}).
\begin{enumerate}
\item A functor $P\colon \BC \xra{} \BA$ with small fibres, that is, such that for each $\BA$-object $A$ the fibre $P^*(A)$ over $A$ is a small set, is a \em discrete fibration over $\BA$ \em if for each $\BA$-morphism $f\colon A'\xra{} PC$ there exists a unique $\BC$-morphism $g\colon C'\xra{}C$ such that $Pg = f$. 

The full subcategory of the slice category $CAT/\BA$ spanned by all discrete fibrations over $\BA$, where $CAT$ stands for the category of all locally small categories, is called the \em category of  discrete fibration over $\BA$ \em and denoted by $DFib(\BA)$. An equivalence in  $CAT/\BA$ is called \em cartesian equivalence\em.

\item Given a discrete fibration $P\colon \BC \xra{} \BA$ its {\em change of base functor} $P^*\colon \mathcal{A}^{\mathsf{op}}\to Set$ is defined by the assignment 
 $(A\xrightarrow{\phi}B)\mapsto ( P^*(B)\xra{P^*(\phi)} P^*(A))$ with as above  $P^*(A)$  the fibre over the $\mathcal{A}$-object $A$ and  $P^*(\phi)$ maps any $b\in P^*(B)$ to  the domain of the unique $\mathcal{C}$-morphism $f$ with codomain $b$ such that $P(f)=\phi$.
\item The \em category $El(P^\ast)$ of elements of  a functor $P^\ast\colon \BA^\op \xra{} \Set$\em, has as 
\begin{itemize}
\item[-] objects all pairs $(A,a)$   with $A\in \ob(\BA)$ and $a\in P^\ast(A)$, 
\item[-] morphisms $(A,a)\xra{\phi}(B,b)$ all $\BA$-morphisms $A\xra{\phi}B$ with $P^\ast(\phi)(b) = a $ 
with composition and identities as in $\BA$,
\item[-] a functor $P\colon El(P^\ast)\xra{} \BA$ given by $\big((A,a)\xra{\phi}(B,b)\big)\mapsto (A\xra{\phi}B)$.
\end{itemize}
\item For any functor $\BA^\op\xra{P^\ast}\Set$ the functor $El(P^\ast)\xra{P} \BA$ is a discrete fibration and the assignment  $P^\ast\mapsto P$ defines an equivalence of the categories $\Set^{\BA^\op}$ and $DFib(\BA)$.
\end{enumerate}

\subsubsection{Kleisli categories}
As shown in the introduction we will be concerned with the endomorphism monoids  of free algebras.
An alternative way   to describe, for any monad $\BBT$ on a category $\BA$, the category of all free algebras $F^\BBT X$ and their homomorphisms is given by the notion of the \em Kleisli category $\BA_\BBT$ \em of $\BBT$  (see e.g. \cite[VI.5]{MacL}). The equivalence of these categories is given by the fact that homomorphisms $F^\BBT X\xra{}F^\BBT Y$  between free algebras  are nothing but the homomorphic extensions of $\BA$-morphisms $X\xra{}F^\BBT Y$. Consequently, the Kleisli category $\BA_\mathbb{T}$ of a  monad $\BBT= (T,\mu,\eta)$ has as
\begin{enumerate}
\item    its objects the  objects of $\BA$  and 
 \item as its hom-sets  $\BA_\mathbb{T}(X,Y)$ the sets $\BA(X,TY)$; its identities $\mathbf{1}_X$ are the $\BA$-morphisms $\eta_X$ and composition  is defined by
$$(Y\xra{ g}TZ)\bullet (X\xra{ f}TY):=X\xra{f }TY\xra{Tg }TTZ\xra{\mu_Z } TZ.$$
\end{enumerate}
The naturality of this concept  in the context of this note is evident: if we would have used in the introduction, instead of $\Set_\sm$, the Kleisli category $\BBT_\sm$  of the free (right) $\sm$-act monad, 
the description of the map $E_X$ would have simply read as $f\mapsto \langle id_X,f \rangle$, while the reconstruction map $L_X$ could simply be described as $(X\xra{\langle id_X ,f \rangle }X\times M)\mapsto  (X\xra{\langle id_X ,f \rangle}X\times M\xra{\pi_M}M)$

\subsection{Two discrete fibrations of an internal category}

From now on one assumes that the finitely complete category $\BC$ is locally small. 
Given   an internal category in $\BC$, that is, a monoid  $\mathsf{M} = (\bar{M},\eta,\mu)$ in the monoidal category $\BC_O$, we denote by  $\BBT_\mathsf{M} = (T_\mathsf{M},\mu_\mathsf{M},\eta_\mathsf{M})$ 
the monad  of its right modules, that is in detail, $T_\mathsf{M}(\bar{A} ) = \bar{A}\ot\bar{M}$,  $(\mu_\mathsf{M})_{\bar{A}}=\bar{A}\ot\bar{M}\ot\bar{M}\xra{id\ot\mu} \bar{A}\ot\bar{M}$ and $(\eta_\mathsf{M})_{\bar{A}} = \bar{A}\simeq \bar{A}\ot\bar{O}\xra{id\ot \eta} \bar{A}\ot\bar{M}$. 
$\widehat{T_\mathsf{M}}$ denotes the composite $T_\mathsf{M}\circ D$, where $D\colon \BC/O\rightarrow\BC_O$  is the embedding.

In view of the introduction  we should be interested in 
\begin{enumerate}
\item   the category  $Conv_\sm = D\downarrow \bar{M}$ of $\bar{M}$-valued morphisms  on $\widehat{\BC}_O$-objects.  (Note that $Conv_\sm$ is locally small as so is $\BC$.) This category is equipped with the forgetful functor $P_\sm\colon Conv_\sm\xra{}\BC/O$ acting as $$\big((\mathsf{f}_A\xra{\alpha}\bar{M}) \xra{\phi}(\mathsf{g}_B\xra{\beta}\bar{M})\big)  \mapsto (A\xra{f}O)\xra{\phi}(B\xra{g}O)$$ (of course, $P_\sm$ has small fibres since $\BC$ is locally small)
\item the endomorphism monoids of free $\BBT_\sm$-modules over $\widehat{\BC}_O$-objects  $\mathsf{f}_A$, i.e., of the monoids
$$\BC_O(\mathsf{f}_A, \mathsf{f}_A\ot \bar{M})=(\BC_O)_{\BBT_\mathsf{M}}(\mathsf{f}_A,\mathsf{f}_A) $$  where $(\BC_O)_{\BBT_\mathsf{M}}$ denotes the Kleisli category of the monad $\BBT_\sm$ and the monoid structure is given by the Kleisli composition $\bullet$ as   $$(\mathsf{f}_A\xra{ \beta}\mathsf{f}_A\ot \bar{M})\bullet (\mathsf{f}_A\xra{\alpha}\mathsf{f}_A\ot \bar{M}) = \mathsf{f}_A\xra{\alpha }\mathsf{f}_A\ot \bar{M}\xra{\beta\ot \bar{M} }\mathsf{f}_A\ot \bar{M}\ot \bar{M}\xra{\mathsf{f}_A\ot \mu} \mathsf{f}_A\ot \bar{M}.$$

These endomorphisms are objects  of the comma category  $D\downarrow \widehat{T_\mathsf{M}}$ and, hence, generate a full  subcategory $End_{\widehat{T_\mathsf{M}}}$ of $D\downarrow \widehat{T_\mathsf{M}}$.  
In more detail   $End_{\widehat{T_\mathsf{M}}}$ has as  objects all 2-cells  $(\mathsf{f}_A\xra{\alpha}\mathsf{f}_A\ot \bar{M}) $ and as morphisms $(\mathsf{f}_A\xra{\alpha}\mathsf{f}_A\ot \bar{M})\rightarrow (\mathsf{g}_B\xra{\beta}\mathsf{g}_B\ot \bar{M})$ pairs of  2-cells $(\mathsf{f}_A\xrightarrow{\sigma}\mathsf{g}_B,\mathsf{f}_A\xrightarrow{\tau}\mathsf{g}_B)$ 
making the following diagram commute. (Note that $End_{\widehat{T_\mathsf{M}}}$ is locally small as so is $\BC$.)
\begin{equation*}
\xymatrix{
\mathsf{f}_A\ar[r]^{\alpha}\ar[d]_{\sigma}& \mathsf{f}_A\ot\bar{M}\ar[d]^{\tau\ot id}\\
\mathsf{g}_B\ar[r]_{\beta}&\mathsf{g}_B\ot\bar{M}
}
\end{equation*}
 $Aut_{\widehat{T_\mathsf{M}}}$ then denotes the full subcategory of  $End_{\widehat{T_\mathsf{M}}}$  spanned by all $\mathsf{f}_A\xra{\alpha}\mathsf{f}_A\ot \bar{M}$  which are   $(\BC_O)_{\BBT_\mathsf{M}}$-automorphisms.  
 \end{enumerate}

 Generalizing the maps $E$ and $L$ of the introduction we will  use the following maps for each $\mathsf{f}_A$ in $\widehat{\BC}_O$  where, obviously, $\bar{\tilde{\alpha}} =\alpha$.
 \begin{enumerate}
\item  $E_{\mathsf{f}_A}\colon\BC_O(\mathsf{f}_A, \bar{M})\xra{}\BC_O(\mathsf{f}_A,\mathsf{f}_A\ot \bar{M} ) $ is the map given by (for the identity see Item 2 of Section \ref{fact:pseudo})
$$(\mathsf{f}_A\xra{ \alpha}\bar{M})\mapsto \tilde{\alpha}:=\mathsf{f}_A\xra{ \langle id, \alpha\rangle}\mathsf{f}_A\ot \bar{M} =\mathsf{f}_A\xra{\Delta}\mathsf{f}_A\ot\mathsf{f}_A\xra{id \ot \alpha}\mathsf{f}_A\ot  \bar{M}.$$
\item  $L_{\mathsf{f}_A}\colon  \BC_O(\mathsf{f}_A,\mathsf{f}_A\ot \bar{M} ) \xra{} \BC_O(\mathsf{f}_A, \bar{M})$ is the map given by ($\bar{\alpha}$ is a 2-cell by Item 4 of Section \ref{fact:pseudo})
$$(\mathsf{f}_A\xra{ \alpha}\mathsf{f}_A\ot \bar{M})\mapsto \bar{\alpha}:=\mathsf{f}_A\xra{ \alpha}\mathsf{f}_A\ot\bar{M}\xra{\pi_d^A}\bar{M}.$$
\end{enumerate}

\subsubsection{The convolution fibration}\label{sec:Mval}

Recall the concept of a convolution monoid as follows: the hom-functor $\BC^\op\times\BC\xra{}\Set$ {of a monoidal (locally small) category $\BC$} is a monoidal functor and, because $(Comon(\BC))^\op\times Mon(\BC) = \Mon(\BC^\op\times\BC)$, it lifts to a functor $(Comon(\BC))^\op\times Mon(\BC) )\xra{Conv}\Mon$.
In other words: if  $(M,M\ot M\xra{m}M,I\xra{e}M)$ a monoid and $(C,C\xra{\delta}C\ot C, C\xra{\epsilon}I)$ a comonoid  in a monoidal category $(\BC,\ot, I)$, then the hom-set $\BC(C,M)$ becomes a monoid  $Conv\big((C,\delta,\epsilon), (M,m,e)\big)$ with unit given as $C\xra{\epsilon}I\xra{e}M$ and with multiplication $\ast$ defined by $(C\xra{f}M)\ast (C\xra{g}M) = C\xra{\delta}C\ot C\xra{f\ot g}M\ot M\xra{m}M$.

Since every object $C$ of a category  with finite products carries a unique comonoid structure whose comultiplication is the diagonal $C\xra{\Delta}C\times C$ and  whose counit is the unique morphism $C\xra{!}1$, one not only sees that the convolution monoids for $\BC=\Set$ are the powers of monoids, but one also obtains  by Lemma \ref{lem:crux} the following

\begin{Fact}\label{fact:}\rm
Given an internal category in $\BC$ considered as a monoid $\mathsf{M} = (\bar{M},\mu,\eta)$ in $\BC_O$ then, for each ${\widehat{\BC}_O}$-objects $\mathsf{f}_A $, the hom-set $\BC_O(\mathsf{f}_A,\bar{M})$ carries the structure of a monoid defined by $\phi\ast\psi := \mathsf{f}_A\xra{\Delta}\mathsf{f}_A\ot \mathsf{f}_A\xra{\phi\ot\psi} \bar{M}\ot \bar{M}\xra{\mu} \bar{M} $ and $e_{\mathsf{f}_A}:= \mathsf{f}_A\xra{f}\bar{O}\xra{\eta}\bar{M}$ and the functor $P^\ast:=
{\widehat{\BC}_O}^{\,\,\,\op}\hookrightarrow\BC_O^{\,\,\op}\xra{\BC_O(-,\bar{M})}\Set$ factors over the category $\Mon$. As one easily sees $El(P^\ast) = Conv_\sm$. This explains the notation chosen for this category: its objects are  the elements of the respective convolution monoids.
\end{Fact}

\begin{Proposition}\label{lem:8} For every internal category $\mathsf{M}$ in $\BC$ and every $\mathsf{f}_A$ in $\BC/O$ the endomorphism monoid $\sm_A(\mathsf{f}_A,\mathsf{f}_A)$  coincides with the convolution  monoid  $ \BC_O(\mathsf{f}_A,\bar{M})$.
\end{Proposition}

\begin{Proof}
By definition of the category   $\sm_A$  (see Section \ref{intcat}) the endomorphism set 
 $\sm_A(\mathsf{f}_A,\mathsf{f}_A) \linebreak = \{A\xra{\alpha}M\mid d\circ\alpha =f = c\circ\alpha\}$ coincides  with $ \BC_O(\mathsf{f}_A,\bar{M} )=   {P_\sm}^\ast(\mathsf{f}_A) $.
  The required equality $\alpha\circ\beta= \alpha\ast\beta$ follows by Item 2 of Section \ref{fact:pseudo}, since it is in detail nothing but\\
$ A\xra{\langle \alpha,\beta\rangle}M\times_OM\xra{\mu} M = |\mathsf{f}_A\xra{\langle \alpha,\beta\rangle}\bar{M}\ot\bar{M}\xra{\mu} \bar{M}| =  %
A\xra{\Delta}A\times_O A\xra{\alpha\times\beta} {M}\times_O {M}\xra{\mu} {M}. $
\end{Proof}

 By Section \ref{discfib} above one also obtains the following results.
\begin{Lemma}\label{lem:7}
The functor   $   {P_\sm}\colon Conv_\sm \xra{}\BC/O$ is a discrete fibration   whose change of base functor  is  ${P_\sm}^\ast:={\widehat{\BC}_O}^{\,\,\,\op}\hookrightarrow\BC_O^{\,\,\op}\xra{\BC_O(-,\bar{M})}\Set$, and this functor factors over the category $\Mon$ of monoids. 
\end{Lemma}

\subsubsection{The free-module-endomorphism fibration}

Recall that for every 2-cell $\mathsf{f}_A\xra{ \alpha}\mathsf{f}_A\ot \bar{M}$ the 2-cells  $\bar{\alpha}=\mathsf{f}_A\xra{ \pi_d^A\circ\alpha}\bar{M}$ and $\alpha':=\mathsf{f}_A\xra{ \pi_f\circ\alpha}\mathsf{f}_A$  form the only pair of 2-cells $(\mathsf{f}_A\xra{\xi}\mathsf{f}_A ,\mathsf{f}_A\xra{  \upsilon}\bar{M})$ \em presenting $\alpha$\em, that is,  such that  $\mathsf{f}_A\xra{ \langle \xi, \upsilon\rangle}\mathsf{f}_A\ot \bar{M}$ exists and equals $\alpha$.

\begin{Definition}\label{def:}\rm
A 2-cell $\mathsf{f}_A\xra{ \alpha}\mathsf{f}_A\ot \bar{M}$ is \em simply presented \em if it is presented by $(id_A, \bar{\alpha})$. $^{sp}End_{\widehat{T_\mathsf{M}}}$ denotes the  full subcategory of $End_{\widehat{T_\mathsf{M}}}$ spanned by   all simply presented    {endomorphisms} while  $^{sp}Aut_{\widehat{T_\mathsf{M}}}$ is its full subcategory spanned by the automorphisms.
\end{Definition}

\begin{Remark}\label{rem:sp}\rm 
A 2-cell  $\mathsf{f}_A\xra{ \alpha}\mathsf{f}_A\ot \bar{M}$ is simply presented iff any of the following holds.
\begin{enumerate}
\item $ \alpha = \mathsf{f}_A\xra{\langle id, \upsilon\rangle} \mathsf{f}_A\ot \bar{M} = \mathsf{f}_A\xra{\Delta} \mathsf{f}_A\ot\mathsf{f}_A \xra{ id\ot \upsilon} \mathsf{f}_A\ot \bar{M}$ for some 2-cell $\mathsf{f}_A\xra{\upsilon}\bar{M}$
\item $\alpha$ belongs to the image    { $\mathcal{E}_{\mathsf{f}_A}$} of the map $E_{\mathsf{f}_A}$.
\item {$E_{\mathsf{f}_A}(L_{\mathsf{f}_A}(\alpha) )=\widetilde{\bar{\alpha}}=\alpha$.}
\end{enumerate}
\end{Remark}

\begin{Lemma}\label{lem:full}
If $(\sigma,\tau)\colon  (\mathsf{f}_A\xra{\langle id, \alpha\rangle} \mathsf{f}_A\ot \bar{M})\xra{}  (\mathsf{g}_B\xra{\langle id, \beta\rangle} \mathsf{g}_B\ot \bar{M}) $ is a morphism in $^{sp}End_{\widehat{T_\mathsf{M}}}$ then $\sigma = \tau$.
\end{Lemma}
\begin{Proof}
 If $(\sigma,\tau)\colon  (\mathsf{f}_A\xrightarrow{\alpha}\mathsf{f}_A\ot \bar{M})\xra{} (\mathsf{g}_B\xrightarrow{\beta}  \mathsf{g}_B \ot\bar{M}$) is  a morphism in $^{sp}End_{\widehat{T_\mathsf{M}}}$, then  $\tau=\tau\circ\pi_f\circ \alpha=\pi_f\circ (id\otimes\tau)\circ \alpha=\pi_g\circ \beta\circ \sigma=\sigma$. 
\end{Proof}

Correspondingly, we may denote morphisms in this category simply by $\sigma$.

\begin{Fact}\label{fact:13}\rm
For each $\mathsf{f}_A\xra{\sigma}\mathsf{g}_B $ in $\BC/O$ there is the map $\mathcal{E}_{\mathsf{g}_B}\xra{\sigma^\ast}\mathcal{E}_{\mathsf{f}_A}$ given by the assignment
 $$(\mathsf{g}_B\xra{\beta}  \mathsf{g}_B\ot\bar{M})\mapsto 
(\mathsf{f}_A\xra{\langle id_A, \pi_d\circ\beta\circ\sigma  \rangle}  \mathsf{f}_A\ot   \bar{M})$$
and the assignment $(\mathsf{f}_A\xra{\sigma}\mathsf{g}_B)\mapsto     {Q_\sm}^\ast(\mathsf{g}_B):=
\mathcal{E}_{\mathsf{g}_B}\xra{\sigma^\ast}\mathcal{E}_{\mathsf{f}_A} =   {Q_\sm}^\ast(\mathsf{f}_A)$ 
defines a functor $(\BC/O)^\op\xra{   {Q_\sm}^\ast}\Set$. 
 
 $El(   {Q_\sm}^\ast)$ has as objects all pairs $(\mathsf{f}_A,\langle id_A,\alpha\rangle )   $ with a 2-cell $\mathsf{f}_A\xra{\alpha}\bar{M}$, while 
a morphism $(\mathsf{f}_A\xra{\langle id,\alpha\rangle}\mathsf{f}_A\ot \bar{M})\xra{\sigma}(\mathsf{g}_B\xra{\langle id,\beta\rangle}\mathsf{g}_B\ot \bar{M})$ in $El(   {Q_\sm}^\ast)$ is a 2-cell $\mathsf{f}_A\xra{\sigma}\mathsf{g}_B $ satisfying the condition $\beta\circ\sigma = \alpha$.
The obvious commutativity of the following diagram shows that this condition is equivalent to the fact that $\sigma$ is a morphism in ${^{sp}End_{\widehat{T_\mathsf{M}}}}$. 

\begin{equation*}\label{diag:5}
\begin{aligned}
\xymatrix@=2.5em{
 \mathsf{f}_A \ar[r]_{\Delta}\ar[d]_{\sigma}\ar@/^2pc/@{->}[rr]_{\langle id,\alpha\rangle}&  \mathsf{f}_A\ot \mathsf{f}_A \ar[d]^{\sigma\ot\sigma} \ar[r]_{id\ot \alpha }&{  \mathsf{f}_A\ot \bar{M}} \ar[d]^{ \sigma\ot id }  \\
\mathsf{g}_B \ar[r]^{\Delta }\ar@/_2pc/@{->}[rr]^{\langle id,\beta\rangle}
 &  \mathsf{g}_B\ot \mathsf{g}_B   \ar[r]^{  id\ot \beta}          & \mathsf{g}_B\ot \bar{M}
}
\end{aligned}
\end{equation*}
Thus,  $El(   {Q_\sm}^\ast)$ is isomorphic to the category ${^{sp}End_{\widehat{T_\mathsf{M}}}}$. 
$   {Q_\sm}$ maps $(\mathsf{f}_A\xra{\langle id,\alpha\rangle}\mathsf{f}_A\ot \bar{M})\xra{\sigma}(\mathsf{g}_B\xra{\langle id,\beta\rangle}\mathsf{g}_B\ot \bar{M})$ to $(A\xra{f}O)\xra{\sigma}  (B\xra{g}O)$. ($Q_\sm$ has small fibres since $\BC$ is locally small.)
\end{Fact}

\begin{Lemma}\label{prop:mon}
The map $E_{\mathsf{f}_A}$ is a homomorphism of monoids from the convolution monoid $\BC_O(\mathsf{f}_A,\bar{M})$ into  the endomorphism  monoid $(\BC_O)_{\BBT_\mathsf{M}}(\mathsf{f}_A,\mathsf{f}_A) $. 
In particular, $\mathcal{E}_{\mathsf{f}_A}$ is a submonoid  of   the endomorphism  monoid $(\BC_O)_{\BBT_\mathsf{M}}(\mathsf{f}_A,\mathsf{f}_A) $. 
\end{Lemma}

\begin{Proof}
First, as easily seen, one has $E_{\mathsf{f}_A}(\mathsf{f}_A\xra{f}\bar{O}\xra{\eta}\bar{M} ) = \mathsf{f}_A\simeq\mathsf{f}_A\ot\bar{O}\xra{id\ot \eta}\mathsf{f}_A\ot\bar{M} $, that is $E_{\mathsf{f}_A}$ preserves  units. $E_{\mathsf{f}_A}$ preserves multiplication iff $  E_{\mathsf{f}_A}( \alpha\ast\beta) = E_{\mathsf{f}_A}( \alpha)\bullet E_{\mathsf{f}_A}(\beta)$,  for any pair  of $\BC_O$-morphisms $\alpha,\beta\colon \mathsf{f}_A\xra{}\bar{M}$, that is, iff the following diagram commutes in $\BC_O$.  But this is  the case by Item 3 of Section \ref{fact:pseudo}.

\begin{equation*}\label{diag:6}
\begin{aligned}
\xymatrix@=2.3em{
\mathsf{f}_A\ar[r]^\Delta\ar[d]_\Delta & \mathsf{f}_A\ot \mathsf{f}_A   \ar[r]^{\Delta  \ot\beta} &{\mathsf{f}_A\ot \mathsf{f}_A\ot \bar{M} } \ar[d]^{\ id\ot\alpha\ot id_M }  \\
\mathsf{f}_A\ot \mathsf{f}_A\ar[d]_{id \ot \Delta}  
&& \mathsf{f}_A\ot  \bar{M}\ot \bar{M}\ar[d]^{id\ot m}\ar@{=}[dl]\\
 \mathsf{f}_A\ot \mathsf{f}_A\ot \mathsf{f}_A\ar[r]_{  id\ot\alpha\ot\beta}   & \mathsf{f}_A\ot \bar{M}\ot \bar{M}\ar[r]_{id\ot m}& \mathsf{f}_A\ot \bar{M}
}
\end{aligned}
\end{equation*}
\end{Proof}

By Section \ref{discfib} we so obtain
\begin{Proposition}\label{prop:}
The functor $   {Q_\sm}\colon {^{sp}End_{\widehat{T_\mathsf{M}}}}\xra{}\BC/O$ is a discrete fibration whose fibres are the sets $\mathcal{E}_{\mathsf{f}_A}$ of simply presented elements and whose change of base functor factors over $\Mon$.
\end{Proposition}

The final and somewhat  surprising result of this section requires the following generalization of Lemma \ref{lem:full}:
\begin{Lemma}\label{lem:lem_corefl}
For each $End_{\widehat{T_{\mathsf{M}}}}$-morphism $(\mathsf{g}_B\xrightarrow{{\beta=\langle id,\bar{\beta}\rangle}}\mathsf{g}_B\otimes\bar{M})\xrightarrow{(\phi,\psi)}(\mathsf{f}_A\xrightarrow{\alpha}\mathsf{f}_A\otimes\bar{M})$ with domain in ${}^{sp}End_{\widehat{T_{\mathsf{M}}}}$ one has $\psi=\phi_\alpha:=\mathsf{g}_B\xrightarrow{\phi}\mathsf{f}_A\xrightarrow{\alpha}\mathsf{f}_A\otimes\bar{M}\xrightarrow{\pi_f}\mathsf{f}_A$.
\end{Lemma}
\begin{Proof}
The claim follows by the following  calculation where the second equality comes from the definition of $End_{\widehat{T_{\mathsf{M}}}}$-morphisms and the third one holds by definition of $\ot$.
 \begin{equation*}
\begin{array}{lll}
\phi_{\alpha}&=&\mathsf{g}_B\xrightarrow{\phi}\mathsf{f}_A\xrightarrow{\alpha}\mathsf{f}_A\otimes\bar{M}\xrightarrow{\pi_f}\mathsf{f}_A=\mathsf{g}_B\xrightarrow{{\langle id,\bar{\beta}\rangle}}\mathsf{g}_B\otimes\bar{M}\xrightarrow{\psi\otimes id}\mathsf{f}_A\otimes\bar{M}\xrightarrow{\pi_f}\mathsf{f}_A\\
&=&\mathsf{g}_B\xrightarrow{{\langle id,\bar{\beta}\rangle}}\mathsf{g}_B\otimes\bar{M}\xrightarrow{\pi_g}\mathsf{g}_B\xrightarrow{\psi}\mathsf{f}_A=\mathsf{g}_B\xrightarrow{id}\mathsf{g}_B\xrightarrow{\psi}\mathsf{f}_A=\psi.
\end{array}
\end{equation*}
\end{Proof}
 \begin{Proposition}\label{prop:corefl}
The category $^{sp}End_{\widehat{T_\mathsf{M}}}$ is a full coreflective subcategory of $End_{\widehat{T_\mathsf{M}}}$.
\end{Proposition}
\begin{Proof}
Observe that  the 2-cell 
$C(\alpha):= \langle id, \pi_d\circ\alpha\rangle\colon \mathsf{f}_A\xra{{\langle id, \bar{\alpha}\rangle}}\mathsf{f}_A\ot \bar{M}$ belongs to $^{sp}End_{\widehat{T_\mathsf{M}}}$, for each $End_{\widehat{T_\mathsf{M}}}$-object $(\mathsf{f}_A\xra{\alpha}\mathsf{f}_A\ot \bar{M})$.

Since by Item 4 (b) of Section \ref{fact:pseudo}   $\mathsf{f}_A \xra{\alpha }\mathsf{f}_A\ot \bar{M}\xra{\pi_f}\mathsf{f}_A$ is a 2-cell, the pair $(id, \pi_f\circ\alpha)$ will be  a $End_{\widehat{T_\mathsf{M}}}$-morphism $C(\alpha)\xra{c_\alpha}\alpha $ if the following diagram commutes.
\begin{equation*}\label{diag:7}
\begin{aligned}
\xymatrix@=3em{
   \mathsf{f}_A  \ar[r]^{{\langle id,\bar{\alpha}\rangle}  } \ar[d]_{ id }   & {\mathsf{f}_A\ot \bar{M}}\ar[d]^{(\pi_f\circ\alpha)\ot id} \\ 
\mathsf{f}_A  \ar[r]_{\alpha}      &   \mathsf{f}_A\ot \bar{M}
}
\end{aligned}
\end{equation*}
But this is easily seen to be true.

We claim that $c_\alpha$ in fact is a coreflection.  We thus have to show that, for any morphism $End_{\widehat{T_\mathsf{M}}}$-morphism 
$(\mathsf{g}_B\xra{{\beta}}\mathsf{g}_B\ot\bar{M})\xra{(\phi,\psi)}(\mathsf{f}_A\xra{{\alpha}}\mathsf{f}_A\ot\bar{M})$ with an $^{sp}End_{\widehat{T_\mathsf{M}}}$-object $(\mathsf{g}_B\xra{{\beta}}\mathsf{g}_B\ot\bar{M})$, there exists a unique $^{sp}End_{\widehat{T_\mathsf{M}}}$-morphism 
$(\mathsf{g}_B\xra{{\beta}}\mathsf{g}_B\ot\bar{M})\xra{(\sigma,\sigma)} C(\alpha)$ such that $c_\alpha\circ (\sigma,\sigma) = (\phi,\psi)$. By the preceding Lemma this condition is equivalent to
$$c_\alpha\circ (\sigma,\sigma)=(id, \pi_f\circ\alpha)\circ (\sigma,\sigma) = (\sigma, \pi_f\circ\alpha\circ\sigma) = (\phi,\phi_\alpha)=  (\phi, \pi_f\circ\alpha\circ\phi)$$
and, hence, satisfied obviously.
\end{Proof}

\subsubsection{Comparing the fibrations}

As seen above we have $   {Q_\sm}^\ast(\mathsf{f}_A)= E_{\mathsf{f}_A}(   {P_\sm}^\ast({\mathsf{f}_A}))$ and (see Remark \ref{rem:sp} above) a 2-cell $\mathsf{f}_A\xra{ \alpha}\mathsf{f}_A\ot \bar{M}$ is simply presented iff $E_{\mathsf{f}_A}(L_{\mathsf{f}_A}(\alpha)) =\widetilde{\bar{\alpha}}=\alpha$.

Moreover, we can restrict the map  $L_{\mathsf{f}_A}$ to a map $   {Q_\sm}^\ast({\mathsf{f}_A})   \xra{}    {P_\sm}^\ast({\mathsf{f}_A})$ and observe that,  for  each 2-cell $\mathsf{f}_A\xra{\alpha }\bar{M}$, one has
 $L_{\mathsf{f}_A}(E_{\mathsf{f}_A}(\alpha)) = \bar{\tilde{\alpha}} =\overline{ \langle id,\alpha\rangle} = \alpha$. 

One so obtains the first of the following statements; the second one follows by a straightforward calculation, while the third one is Lemma \ref{prop:mon}.  The final item then follows from Lemma \ref{lem:7} by Items 2 and 3.

\begin{Lemma}\label{lem:comparison}
\begin{enumerate}
\item $L_{\mathsf{f}_A}$ is the inverse of $E_{\mathsf{f}_A}$.
\item  The family $(E_{\mathsf{f}_A})_{\mathsf{f}_A}\colon    {P_\sm}^\ast\Rightarrow    {Q_\sm}^\ast$ is a natural isomorphism with $(L_{\mathsf{f}_A})_{\mathsf{f}_A}$ as  its inverse. 
\item Each map $E_{\mathsf{f}_A}$ is an isomorphism of monoids $   {P_\sm}^\ast(\mathsf{f}_A)  \xra{}    {Q_\sm}^\ast (\mathsf{f}_A) $ with $L_{\mathsf{f}_A}$ as its inverse.
\item  For  each morphism $\mathsf{f}_A\xra{\sigma}\mathsf{g}_B $ in $\BC/O$ the map $   {Q_\sm}^\ast(\sigma)$ is a homomorphism of monoids $   {Q_\sm}^\ast(\mathsf{g}_B)\xra{}   {Q_\sm}^\ast(\mathsf{f}_A)$.
\end{enumerate}
\end{Lemma}

In view of Section \ref{discfib}  we  can summarize the above alternatively as follows.
\begin{Proposition}\label{prop:15}
For any internal category $\mathsf{M}= (O,M,d,c,\eta,\mu)$  in a finitely complete category $\BC$ the following hold.
\begin{enumerate}
\item
The functors $   {P_\sm}\colon Conv_\sm \xra{}\BC/O$ and $   {Q_\sm}\colon {^{sp}End_{\widehat{T_\mathsf{M}}}} \xra{}\BC/O$ are discrete fibrations whose change of base functors $   {P_\sm}^\ast$ and $   {Q_\sm}^\ast$ are naturally isomorphic and factor over the category $\Mon$ of monoids.
\item 
The respective cartesian isomorphism $E\colon  Conv_\sm \xra{} {^{sp}End_{\widehat{T_\mathsf{M}}}}$ and its inverse \linebreak $L\colon  {^{sp}End_{\widehat{T_\mathsf{M}}}}\xra{} Conv_\sm$ act as follows:
\begin{eqnarray*}
E\big((\mathsf{f}_A\xra{\alpha}\bar{M})\xra{\phi}(\mathsf{g}_B\xra{\beta}\bar{M}) \big) & = & (\mathsf{g}_B\xra{\tilde{\beta}}\mathsf{g}_B\ot \bar{M}
)\xra{\phi}(\mathsf{f}_A\xra{\tilde{\alpha}}\mathsf{f}_A\ot \bar{M})\\
L\big((\mathsf{f}_A\xra{\alpha}\mathsf{f}_A\ot \bar{M} )\xra{(\phi,\phi)}(\mathsf{g}_B\xra{\beta}\mathsf{g}_B\ot \bar{M}
)\big) 
& = & (\mathsf{f}_A\xra{\bar{\alpha}}\bar{M})\xra{\phi}(\mathsf{g}_B\xra{\bar{\beta}}\bar{M}) 
\end{eqnarray*}
\end{enumerate}
\end{Proposition}

\begin{Remark}\label{rem:21}\rm
By Proposition \ref{lem:8} the convolution monoids $   {P_\sm}^\ast(\mathsf{f}_A)$ are groups if the categories $\sm_A$ are groupoids, and this clearly is the case if the internal category under consideration is an internal groupoid (see Fact \ref{prop:onetoone}). We so obtain: 
if $ (O,M,d,c,\eta,\mu),\iota)$  is an internal groupoid in $\BC$ then the functors $   {P_\sm}^\ast$ and $   {Q_\sm}^\ast$ even factor over the category $Grp$ of groups.
\end{Remark}

With $E$ and $L$  the isomorphisms of  Proposition \ref{prop:15} and $C$  the coreflector of Proposition \ref{prop:corefl} one obtains
\begin{Corollary}\label{corr:}
Let $\mathsf{M} = (O,M,d,c,\eta,\mu)$ be an internal category in the finitely complete category $\BC$. Then
the following hold:
\begin{enumerate}
\item The functor $\hat{E} =  Conv_\sm \xra{E }{^{sp}End_{\widehat{T_\mathsf{M}}}}\hookrightarrow End_{\widehat{T_\mathsf{M}}} $ is left adjoint to the functor $\hat{L}\colon End_{\widehat{T_\mathsf{M}}}\xra{C }{^{sp}End_{\widehat{T_\mathsf{M}}}}\xra{L} Conv_\sm$.
\item If $(\mathsf{M},\iota)$ even is an internal groupoid in $\BC$ then this adjunction induces, an adjunction $E'\dashv L'$ with $L' = Aut_{\widehat{T_\mathsf{M}}}\hookrightarrow End_{\widehat{T_\mathsf{M}}}\xra{L } Conv_\sm$, the restriction of the functor  $\hat{L}$, and  $E' \colon Conv_\sm\xra{} Aut_{\widehat{T_\mathsf{M}}}$ the corestriction of $\hat{E}$.
\end{enumerate}
\end{Corollary}

\subsection{Functoriality}

We finally investigate to what extent the assignments $\sm\mapsto P_\sm$ and $\sm\mapsto Q_\sm$ are functorial. For doing so it seems appropriate to define  the category ${I\!ntCat}$  \em of all internal categories \em in analogy to the category of all modules and also  the \em category  of all discrete fibrations \em ${DFib}$. We will  need the observation that the construction of Remark \ref{rem:Catin}  defines a functor $Cat_{in}\colon {Lex}\to CAT$ where ${Lex}$ denotes the category of all finitely complete  {locally small categories and finite limit preserving functors and \em CAT \em that of all locally small categories}. 

Since these constructions involve categories of categories let us clarify our conventions concerning the sizes of the categories that are considered below: we will here use the term ``category" also for structures which are roughly like ``metacategories" in \cite{MacL} or could be called  ``very large categories'', such as $Lex$, $CAT$, ${I\!ntCat}$ or $DFib$ (see below). In particular such categories are not assumed locally small.\footnote{More formally, for the sake of the reader interested in the foundations of the theory of categories: let us choose e.g. as our foundational system ZFC $+$ Grothendieck universes (that is, for every set there exists a Grothendieck universe containing it), and  let us called ``category'' any set-theoretic model of the axioms of categories. Let $U_1\in U_2\in U_3$ be Grothendieck universes. Let us fix the terminology, as already but  informally used in this note: a set which belongs to $U_1$ (resp. $U_2$, resp. $U_3$) is referred to as {\em small} (resp. {\em large}, resp. {\em very large}). A {\em small} (resp. {\em large}, resp. {\em very large}) category is then a category whose sets of objects and arrows are small (resp. large, resp. very large). A large category is {\em locally small} when its hom-sets are small. E.g. the categories of small sets $Set$ and of small categories $Cat$ are locally small. The categories $Lex$, $CAT$, ${I\!ntCat}$ and $DFib$  are very large, while  only ``locally large'', that is,   their hom-sets are  large.}

\begin{description}
\item[The category ${I\!ntCat}$ of internal categories] 
 has as objects the pairs $(\mathcal{C},\mathsf{M})$ formed by a finitely complete  category $\mathcal{C}$ and  an internal category $\mathsf{M }$ in $\mathcal{C} $, and as morphisms the pairs $(\mathcal{C},\mathsf{M})\xrightarrow{(K,(F_o,F_m))}(\mathcal{C}',\mathsf{M}')$ consisting of a finite limit  preserving functor $K\in {Lex}(\mathcal{C},\mathcal{C}')$ and an internal functor $ Cat_{in}(K)(\mathsf{M})\xrightarrow{(F_o,F_m)}\mathsf{M}'$ (see Remark \ref{rem:Catin}). In other words, ${I\!ntCat}$  is the category obtained from the functor {$Cat_{in}\colon {Lex}\to CAT $ } using the Grothendieck construction (see e.g.~\cite[Theorem~1.10.7, p.~111]{Jacobs}). 

Given a morphism $(\mathcal{C},\mathsf{M})\xrightarrow{(K,(F_o,F_m))}(\mathcal{C}',\mathsf{M}')$  in ${I\!ntCat}$ we denote by $K/F_o\colon \BC/O\to\BC'/O'$ the functor mapping $\big((C\xra{f}O)\xra{\alpha}(D\xra{g}O)\big)$ to $\big(( K(C)\xra{K(f)}K(O)\xra{F_o}O')\xra{K(\alpha)}(K(D) \xra{K(g)}K(O)\xra{F_o}O'  \big)$.

\item[The category ${DFib}$ of discrete fibrations] is the full subcategory of the category of arrows in {$CAT$}  spanned by all discrete fibrations; in other words,    if $\mathcal{A}\xra{P} \mathcal{C}$ and  $\mathcal{B}\xra{Q} \mathcal{D}$ are discrete fibrations, then a morphism $P\to Q$ is a pair of functors $(\mathcal{A}\xra{F}\mathcal{B},\mathcal{C}\xra{G}\mathcal{D})$ such that $Q\circ F = G\circ P$.

For any pair $(\BC,\sm)$ in ${I\!ntCat}$ we denote by
\begin{enumerate}
\item ${Conv}$  the  full subcategory of ${DFib}$ spanned by the discrete fibrations \linebreak $Conv_\sm \xra{P_\sm}\BC/O$ 
\item $^{sp}{End}$ be the  full subcategory of ${DFib}$ spanned by the discrete fibrations  $^{sp}End_{\widehat{T_{\mathsf{M}}}}\xra{Q_\sm}\BC/O$.
\end{enumerate}

\end{description}
One now can define functors ${I\!ntCat}\xra{}{DFib}$ as follows:
\begin{description}
\item[P\ \ \ ] For any 
$(\mathcal{C},\mathsf{M})\xrightarrow{(K,(F_o,F_m))}(\mathcal{C}',\mathsf{M}')$  in ${I\!ntCat}$
let  $$[K,(F_o,F_m)]\colon Conv_\sm\to  Conv_{\sm'}$$ be the functor   acting
 \begin{itemize}
\item  on objects by $(\mathsf{f}_A\xrightarrow{\alpha}\bar{M})\mapsto (\mathsf{(F_o\circ K(f))}_{K(A)}\xrightarrow{F_m\circ K(\alpha)}\bar{M}')$
\item and on morphisms by $\sigma\mapsto K(\sigma)$.
\end{itemize}

{It is now easy to see 
that the following diagram commutes and, hence, the pair $( [K,(F_o,F_m)]  , K/F_o )$ is a morphism $P_\sm\xra{}P_{\sm'}$ in  ${Conv}  $.\begin{equation*}\label{diag:8}
\begin{aligned}
\xymatrix@=2.5em{
   Conv_\sm \ar[rr]^{\!\!\! [K,(F_o,F_m)]  }\ar[d]_{P_\sm  }&&{ 
     Conv_{\sm'}} \ar[d]^{P_{\sm'} }  \\
  \BC/O  \ar[rr]_{ K/F_o}         & &    \BC'/O'
}
\end{aligned}
\end{equation*}
In fact, commutativity of this diagram is equivalent to saying that for any morphism 
 $(\mathsf{f}_A\xrightarrow{\alpha}\bar{M})\xrightarrow{\sigma}(\mathsf{g}_B\xrightarrow{\beta}\bar{M})$ in $Conv_{\mathsf{M}}$, i.e., for any commutative diagram
\begin{equation*}
\xymatrix{
&A  \ar[d]^{\sigma}\ar[rd]^{f}\ar[ld]_{f}&\\
O &\ar[l]_{g}\ar[r]^{g}B\ar[d]^{\beta}&O\\
&\ar[lu]^{c}\ar[ru]_{d}\ar[lu]^{c}\ar[ru]_{d}M&
}
\end{equation*}
with $\alpha= \beta\circ\sigma$ the following diagram commutes
\begin{equation*}
\xymatrix{
&\ar[d]^{K(\sigma)}\ar[ld]_{K(f)}K(A)\ar[rd]^{K(f)}&\\
K(O)\ar[d]_{F_o}&\ar[l]_{K(g)}\ar[r]^{K(g)}K(B)\ar[d]^{K(\beta)}&K(O)\ar[d]^{F_o}\\
O'&\ar[lu]^{K(c)}K(M)\ar[d]^{F_m}\ar[ur]_{K(d)}&O'\\
&\ar[lu]^{c'}\ar[ru]_{d'}M'&
}
\end{equation*}
But this is obvious.}

A further straightforward calculation shows that the assignment 
$$\big((\BC,\sm)\xra{(K,(F_o,F_m))}(\BC',\sm')\big)\to \big(P_\sm\xra{( [K,(F_o,F_m)]  , K/F_o) }P_{\sm'}\big)$$ 
defines a functor 
${I\!ntCat}\xra{P}{Conv}\hookrightarrow {DFib}$.
\end{description}

\begin{description}
\item[Q\ \ \ ] For any $(\mathcal{C},\mathsf{M})\xrightarrow{(K,(F_o,F_m))}(\mathcal{C}',\mathsf{M}')$ 
 in ${I\!ntCat}$ 
 let $$\{K,(F_o,F_m)\}\colon {}^{sp}End_{\widehat{T_{\mathsf{M}}}}\to {}^{sp}End_{\widehat{T_{\mathsf{M}'}}}$$ be the functor  acting 
\begin{itemize}
\item on  objects by  $$((\mathsf{f}_A\xrightarrow{\langle id,\alpha\rangle}\mathsf{f}_A\otimes\bar{M}) \mapsto \big((\mathsf{F_o\circ K(f)})_{K(A)}\xrightarrow{\langle id, F_m\circ K(\alpha)\rangle}(\mathsf{F_o\circ K(f)})_{K(A)}\otimes\bar{M}'\big)$$
\item and on morphisms by $\sigma\mapsto K(\sigma)$.
\end{itemize}

Again by a simple calculation one obtains that the following diagram commutes and, hence, the pair $( \{K,(F_o,F_m)\} , K/F_o )$ is a morphism $Q_\sm\xra{}Q_{\sm'}$ in  $^{sp}{End}$,  
\begin{equation*}\label{diag:9}
\begin{aligned}
\xymatrix@=4em{
  {}^{sp}End_{\widehat{T_{\mathsf{M}}}}   \ar[r]^{\!\!\! \{K,(F_o,F_m)\} }\ar[d]_{Q_\sm  }&{ {}^{sp}End_{\widehat{T_{\mathsf{M}'}}}} \ar[d]^{Q_{\sm'}}\\
  \BC/O  \ar[r]_{ K/F_o}          &    \BC'/O'
}
\end{aligned}
\end{equation*}  

{while a } further  calculation shows that the assignment   
$$\big((\BC,\sm)\xra{(K,(F_o,F_m))}(\BC',\sm')\big) \to \big(Q_\sm\xra{ ( \{K,(F_o,F_m)\}  , K/F_o )  } Q_{\sm'} \big)$$ 
defines a functor ${I\!ntCat}\xra{Q} {}^{sp}{End}\hookrightarrow {DFib}$.

\item[Conclusion] The functors $P$ and $Q$, considered as functors ${I\!ntCat}\xra{}{DFib}$, are naturally isomorphic.
\end{description}

\section{Possible applications}\label{Sec:applications}

\subsection{A ``quantum'' Feistel scheme}

The purpose of this section is a generalization of
 the cryptographic construction of Feistel schemes referred to in the introduction, which is suitable 
 for quantum computing. We therefore briefly recall the original construction as follows,  where plaintexts and ciphertexts are members of $2^{2m}$:
choose,  for some $N\in\N$ (called \em number of rounds\em ) ``key"-maps $f_i\colon 2^m\xra{} 2^m$ for each $1\leq i\leq N$ and form the maps $E_{m}(f_i)\circ\sigma\colon  2^m\times 2^m\xra{} 2^m\times 2^m$, where 
$\sigma\colon (\vec{x},\vec{y})\in 2^{m}\times 2^m\mapsto (\vec{y},\vec{x})\in 2^m\times 2^m$ and  their composition $E:= (E_{m}(f_N)\circ\sigma) \circ \cdots \circ  (E_{m}(f_1)\circ\sigma)$. Then   the encoded ciphertext  $c_t$ corresponding to some plaintext $t$ is computed   as $E(t)$.

Making this construction implementable  on a  quantum computer requires to find a category $\BC$ with finite products such that, for any internal group $(H,\mu,\eta,\iota)$ in $\BC$ and each morphism $H\xrightarrow{f}H$ in $\mathcal{C}$, the morphism $E_H(f)\colon H\times H\to H\times H$ generalizing the fundamental step of constructing Feistel ciphers is a unitary operator on a finite dimensional Hilbert space (see e.g.~\cite{Bera}).

In order to define a suitable category $\mathcal{C}$ recall first the following facts about the category 
${H\!ilb\!_{fd}}$ of finite dimensional complex Hilbert spaces and linear maps (these are automatically bounded): (a) $(H\!ilb\!_{fd},\otimes,\mathbb{C},\sigma)$ with $\otimes$ the usual tensor product and $\sigma$ the usual twist is a symmetric monoidal category; (b) the usual Hilbert adjoint functor ${}^{\dagger}\colon {H\!ilb\!_{fd}}^\op\to H\!ilb\!_{fd}$ makes it a dagger compact closed category (see e.g.~\cite{Selinger}); (c) $H\!ilb\!_{fd}$ has finite biproducts. By $|-|$ we denote its underlying space functor into the category of complex vector spaces. 

As for any monoidal category we, hence, can consider its category ${}_{coc}{Comon}(H\!ilb\!_{fd})$ of cocommutative comonoids. Its objects are triples $({H},\delta,\epsilon)$ where ${H}$ is a finite dimensional Hilbert space, and where $(|H|,|H|\xrightarrow{\delta}|H|\otimes_\mathbb{C}|H|,|H|\xrightarrow{\epsilon}\mathbb{C})$ is an ordinary  cocommutative (complex) coalgebra with a counit. As is well-known the tensor product of $H\!ilb\!_{fd}$ lifts to ${}_{coc}{Comon}(H\!ilb\!_{fd})$ and here is the binary product; in particular, this category has all finite products (the Hilbert space $\mathbb{C}$ with its trivial coalgebra structure is a terminal object)~\cite{Fox1976}.
A cocommutative comonoid $({H},\delta,\epsilon)$ is a {\em $\dagger$-Frobenius coalgebra} when $(\delta^{\dagger}\otimes id)(\delta\otimes id)=\delta\circ\delta^{\dagger}=(id\otimes \delta^{\dagger})(\delta\otimes id)$. It is called {\em special} when $\delta$ is an isometry, that is, $\delta^{\dagger}\circ\delta=id$. One immediately sees that $(H,\delta^{\dagger},\epsilon^{\dagger})$ then is a special commutative $\dagger$-Frobenius algebra in $H\!ilb\!_{fd}$  (see~\cite{Coecke}).

 The full subcategory  $\BC_{Hilb}$ of ${}_{coc}{Comon}(H\!ilb\!_{fd})$ generated by the special $\dagger$-Frobenius coalgebras turns out to be suitable for our purpose. It first is closed under binary products and contains $\mathbb{C}$  (a direct consequence of~\cite[Proposition~2.2.6, p.~18]{Abrams} for instance) and, thus, has all finite products.
Moreover, an internal group  in this category is  a tuple $((H,\delta,\epsilon),\mu,\eta,S)$ where $(H,\delta,\epsilon)$ is a special cocommutative $\dagger$-Frobenius coalgebra and $(|H|,\delta,\epsilon,\mu,\eta,S)$ is  an ordinary (finite dimensional) Hopf algebra (see e.g.~\cite[p.~377]{Dasca} concerning the relation between cocommutative Hopf algebras and group objects in the category of cocommutative coalgebras). 
In addition one observes the following:
the set $G({H},\delta,\epsilon)$ of group-like elements, that is, the elements $x\in H$ such that $\delta(x)=x\otimes x$ and $\epsilon(x)=1$,  of a special cocommutative $\dagger$-Frobenius coalgebra $({H},\delta,\epsilon)$  forms an orthonormal basis of the Hilbert space ${H}$ (see~\cite[Corollary~4.7]{Coecke} with~\cite[Propositions~14 and~15]{Abramsky}), and for a morphism of coalgebras $({H},\delta,\epsilon)\xrightarrow{f}({H}',\delta',\epsilon')$ between special $\dagger$-Frobenius coalgebras, $f(G({H},\delta,\epsilon))\subseteq G({H}',\delta',\epsilon')$.

Let now  $(({H},\delta,\epsilon),\mu,\eta,S)$ be an internal group in $\BC_{Hilb}$ and let $f\colon \mathsf{H}\to\mathsf{H}$, with $\mathsf{H}:= ({H},\delta,\epsilon)$, be an endomorphism. Then, using the notations from the Introduction, $E_{\mathsf{H}}(f)\colon \mathsf{H}\otimes\mathsf{H}\to \mathsf{H}\otimes\mathsf{H}$, given in details by $u\otimes v\mapsto \sum_{(u)}u_1\otimes \mu(f(u_2)\otimes v)$ (with $\delta(u)=\sum_{(u)}u_1\otimes u_2$), is an automorphism. In particular, since $E_{\mathsf{H}}(f)(G({H}\otimes {H}))=G({H}\otimes{H})$, $E_{\mathsf{H}}(f)$ is in fact  a unitary operator on ${H}\otimes{H}$ as aspired. 

 Admittedly, we don't know whether the execution of these generalized `ciphers' in the context of quantum computing is of great importance.

\subsection{Higher-order Boolean functions}

Let $\bigsqcup_{x\in X}(G_x,*_x,1_x)$ be the coproduct, in the category of groupoids, of an $X$-indexed family $(G_x,*_x,1_x)_{x\in X}$ of groups, each of them being considered as a one-object groupoid. The domain $d$ and codomain $c$ of this groupoid are both equal to  the canonical projection $\pi\colon \bigsqcup_{x\in X}G_x\to X$, $(x,g)\mapsto x$. Let $\alpha\colon X\to \bigsqcup_{x\in X}G_x$ be a section of $\pi$, that is, $\alpha(x)=(x,\alpha^g(x))$ where $\alpha^g(x)\in G_x$, $x\in X$. 
By Remarks \ref{prop:onetoone} and  \ref{rem:21} the map
$E(\alpha)\colon  \bigsqcup_{x\in X}G_x\to \bigsqcup_{x\in X}G_x$ given by $E(\alpha)(x,g)=(x,\alpha^g(x)*_x g)$ is a permutation. 

In the case where $G_x=\Set(2^{m_x},2^{n_x})$ is the group  {${\Z_2}^{n_x2^{m_x}}$},
 $\alpha\colon X\to \bigsqcup_{x\in X} \Set(2^{m_x},2^{n_x})$ is, in the language of functional programming,   a higher-order Boolean function. Hence our construction  turns such a higher-order Boolean function into a higher-order Boolean permutation. This construction might be of use in functional programming.

\section*{Acknowledgements}


The first author thanks Prof. Sami Harari (Toulon) for fruitful discussions about the Feistel scheme which were the germ of the work presented here.


\begin{thebibliography}{99}

\bibitem{Abrams}
L.E. Abrams,  Frobenius algebra structures in topological quantum field theory and quantum cohomology. PhD Thesis, The Johns Hopkins University (1997).

\bibitem{Abramsky}
S. Abramsky and C. Heunen, H$^*$-algebras and nonunital Frobenius algebras:
first steps in infinite-dimensional categorical quantum mechanics, {\em Clifford Lectures} {\bf 71} (2012), 1--24.

\bibitem{AHS} J. Ad\'amek, H. Herrlich H. and G.E. Strecker,  {\em
Abstract and Concrete Categories}, John Wiley, New York (1990).

\bibitem{Ben}  J. B\'enabou, Introduction to bicategories, in: \em Reports of the Midwest Category Seminar\em. Lecture Notes in Mathematics {\bf 47},  Springer, Berlin, Heidelberg, 1967. p. 1--77.

\bibitem{Bera}
R.K. Bera, {\em The Amazing World of Quantum Computing}, Springer, Singapore (2020).

\bibitem{Brown} R. Brown, From groups to groupoids: a brief survey,  {\em Bull. London Math. Soc.} {\bf 19} (1987), 113--134.

\bibitem{Coecke}
B. Coecke, D. Pavlovic and J. Vicary, A new description of orthogonal bases, 
{\em Math. Structures Comput. Sci.} {\bf 23(3)} (2013), 555--567.

\bibitem{Dasca}
S. D\u{a}sc\u{a}lescu, C. N\u{a}st\u{a}sescu and  \c{S}. Raianu, {\em Hopf Algebras}, Monographs and Textbooks in Pure and Applied Mathematics {\bf 401}, Marcel Dekker (2001).


\bibitem{DPP}  R. Dawson, R. Par\'e and 
 D. Pronk, The span construction, \em Theory and Applications of Categories\em, {\bf 24} (2010), 302--377.
 
 \bibitem{Fox1976}
T. Fox, Coalgebras and cartesian categories, {\em Comm. Alg.} {\bf 4(7)} (1976), 665--667.

\bibitem{Illusie} L. Illusie, \em Complexe Cotangent et D\'eformations II\em, LNM 283, Springer (1972).

\bibitem{Jacobs} B. Jacobs,  \em  Categorical Logic and Type Theory\em, Studies in logic and the foundations of mathematics {\bf 141}, Elsevier (1999).

\bibitem{Riehl} F. Loregian and E. Riehl, Categorical notions of fibration, \em Expositiones Mathematicae \em {\bf 8} 
(2020), 496--514.

%
\bibitem{MacL} S. Mac Lane, {\em Categories for the Working Mathematician}, 2nd ed., Springer, New York (1998).
\bibitem{Menezes}
A.J. Menezes, P.C. van Oorschot  and S.A. Vanstone, {\em Handbook of Applied Cryptography (Fifth ed.)}, CRC Press (2001).

\bibitem{Par} B. Pareigis, Non-additive ring and module theory I., {\em Publ. Math.} {\bf 24}, Debrecen, 189--204  (1977).

\bibitem{Selinger}
P. Selinger, Dagger compact closed categories and completely positive maps. {\em Electronic Notes in Theoretical computer science} {\bf 170} (2007), 139--163.

\bibitem{Street} R. Street, The formal theory of monads, {\em J. Pure Appl. Algebra.} {\bf{2}}  (1972), 149--168.

\bibitem{Toffoli}
T. Toffoli, Reversible computing. In : {\em International colloquium on automata, languages, and programming}, Lecture Notes in Computer Science {\bf 85}. Springer, Berlin, Heidelberg, 1980. p. 632--644.

 \end{thebibliography}
\end{document}